\documentclass[10pt,a4paper]{article}

\usepackage[utf8]{inputenc}
\usepackage{amsmath}
\usepackage{amsfonts}
\usepackage{amssymb}
\usepackage{enumitem}
\usepackage{wrapfig}
\usepackage[headsep=10mm]{geometry}
\usepackage{marginnote}
\usepackage{caption}
\usepackage{subcaption}
\usepackage[pdftex]{graphicx}
\usepackage{color}
\usepackage{lipsum}
\usepackage{bm}
\usepackage{mathrsfs}
\usepackage{afterpage}
\usepackage{footmisc}
\usepackage[hidelinks]{hyperref}
\usepackage{fancyhdr}

\numberwithin{equation}{section}

\DeclareMathOperator{\diff}{d}
\def\MM#1{\boldsymbol{#1}}

\newcommand{\grad}[1]{{\bm{\nabla}{#1}}}

\newcommand{\curl}[1]{\bm{\nabla}\times\bm{#1}}
\newcommand{\lap}[1]{\nabla^2{#1}}
\newcommand{\pfrac}[2]{\frac{\partial{#1}}{\partial{#2}}}
\newcommand{\dx}[1]{\hspace{0.75mm}\mathrm{d}{#1}}

\newcommand{\curlybrackets}[1]{\left\lbrace{#1}\right\rbrace}

\newcommand{\new}[1]{{#1}}

\begin{document} 
\title{Statistical properties of an enstrophy conserving finite element discretisation for the stochastic quasi-geostrophic equation}
\author{Thomas M. Bendall and Colin J. Cotter}
\maketitle

\begin{abstract}
\noindent
A framework of variational principles for stochastic fluid
dynamics was presented by Holm (2015), and these stochastic
equations were also derived by
Cotter et al. (2017).
We present a conforming finite element discretisation for the 
stochastic quasi-geostrophic equation that was derived from this framework.
The discretisation preserves the first two moments of potential vorticity, i.e. the mean potential vorticity and the enstrophy.
Following the work of Dubinkina and Frank (2007), who investigated the statistical mechanics
of discretisations of the deterministic quasi-geostrophic equation, we investigate the statistical
mechanics of our discretisation of the stochastic quasi-geostrophic equation.
We compare the statistical properties of our discretisation with the Gibbs distribution under assumption of these conserved quantities, finding that there is agreement between the statistics under a wide range of set-ups.
\end{abstract}


\section{Introduction}
In recent years, there has been increased interest in
stochastic models of geophysical fluids, as they provide a way
to allow models to \new{express some of the spread in uncertainty of unresolved processes}.
One approach is the variational principles framework introduced by  \cite{holm2015variational} for stochastic fluid dynamics.
In this framework, stochastic perturbations (representing the effect of unresolved motions on the resolved scales) are introduced as perturbations to the velocity field that transports Lagrangian fluid particles.
This is expressed as in the following formula describing stochastic Lagrangian transport,
\begin{equation}
  \mathrm{d}\bm{x}_t =\bm{u}(\bm{x},t)\dx{t}-\sum_i\bm{\Xi}_i(\bm{x})\circ\mathrm{d}W_i,
  \label{eq:dx}
\end{equation}
where $\bm{u}$ is the deterministic velocity field, $\bm{\Xi}_i(\bm{x})$ are time-independent basis functions which determine the spatial correlations in the stochastic component of the velocity field, and $\circ\dx{W_i}$ denotes Stratonovich noise.
The stochastic variational principle leads to Eulerian fluid models, with stochastic multiplicative noise in the form of additional transport terms, that have a conserved potential vorticity; the addition of noise breaks the time-translation symmetry and hence they do not conserve energy.
Unlike other approaches, the introduction of stochasticity through
this framework still \new{has a Kelvin circulation theorem since it preserves the particle relabelling symmetry}. 
This framework can be used to develop potential vorticity conserving stochastic versions of all of the geophysical fluid dynamics models that have Hamiltonian structure: shallow water equations, Boussinesq equations, anelastic, pseudo-compressible, compressible Euler, and so on.
The dynamics resulting from this approach has also 
been derived by \cite{cotter2017stochastic}, who decomposed the
deterministic Lagrangian flow map into slow, large-scale motions
and fast, small-scale fluctuations, before using homogenization theory.
\\
\\
Two important potential applications of this approach are 
modelling unresolved subgrid backscatter onto resolved scales
and in ensemble prediction and ensemble data assimilation, where the stochastic forcing can be used to reach other nearby trajectories that are physically realistic. 
\\
\\
In this work our goal is to study the behaviour of the
stochastic equation \new{derived in \cite{holm2015variational}}.
We will consider one of the simplest members of this stochastic model hierarchy, the stochastic quasi-geostrophic (QG) equation,
\begin{equation}
\label{eq:sqg}
  \dx{q} + \mathrm{d}\bm{x}_t\bm{\cdot}\grad{} q=0,
\quad q:=\nabla^2\psi - \mathcal{F}\psi + f,
\end{equation}
where $q(\bm{x},t)$ is the scalar potential vorticity, $\mathcal{F}$ is the non-dimensional Froude number and $f$ is the non-dimensional Coriolis parameter. 
The stream function $\psi$ is related to the velocity $\bm{u}$ through 
\begin{equation*}
\bm{u}=\grad{}^\perp \psi,
\end{equation*}
where the $\grad{}^\perp$ operator is defined by $\grad{}^\perp:=(-\partial_y, \partial_x)$. \\
\\
Our approach is to study the new stochastic QG equation via numerical
models, thereby proposing a first discretisation to this equation.
In this paper we develop a conforming finite element discretisation for the stochastic QG equation, showing that it conserves total PV and enstrophy.
Having presented a new discretisation, we then study its statistical mechanics in numerical experiments that compare the long time statistics of the model with predictions from an appropriate Gibbs distribution,
which helps us to learn about the properties of the underlying
stochastic QG equation.
The novel results of this work are the presentation
of a discretisation to the stochastic QG equation,
\new{using a standard finite element method with implicit midpoint rule time stepping,} 
and the subsequent investigation of its statistical properties.
Through this approach, we find that \new{the equilibrium distribution in phase space of the numerical model is approximated by a Gibbs distribution.}
We also provide a derivation of the stochastic QG model 
directly from the stochastic variational principle using the Lagrangian
of \cite{holm1998hamilton}, contrasting with the Hamiltonian approach
used in \cite{holm2015variational}.
\\
\\
The rest of this paper is structured as follows. 
In Section \ref{sec:sqg} we derive the stochastic QG equation. 
In Section \ref{sec:discrete QG} we describe the finite element discretisation and derive conservation properties.
In Section \ref{sec:statmech} we review the properties of the Gibbs distribution under the assumption of conservation of
enstrophy, and in Section \ref{sec:numerics} we compare statistics from the numerical discretisation with corresponding statistics from the Gibbs distribution.
Finally we provide a summary and outlook in Section \ref{sec:outlook}.
\section{Stochastic quasi-geostrophic model derivation}
\label{sec:sqg}
In this section we derive the stochastic QG equation from a stochastic variational principle, by adapting the variational principle for QG of \cite{holm1998hamilton} to the stochastic framework of \cite{holm2015variational}. 
This equation was previously only derived by \cite{holm2015variational} using a Poisson bracket approach.
\cite{holm1998hamilton} considered a fluid with velocity $\MM{u}$ 
and density $\rho$, together with a Lagrange multiplier $P$ 
(which we will later use to enforce constant density 
$\rho=\rho_0$ and hence incompressibility).
Following the framework of \cite{holm1998hamilton}, we construct the following action
\begin{align}
\begin{split}
  S[u,\rho,\phi,\MM{p},\MM{q},P] &=
  \int  \ell(\MM{u},\rho,P) \dx{t} + \int_\Omega\phi\left(\mathrm{d}\rho + 
  \left[\MM{u}\dx{t} + \sum_i\bm{\Xi}_i(\bm{x})\circ\mathrm{d}W_i\right] \bm{\cdot}\grad{\rho}\right)
 \\
& \qquad  + \MM{p}\bm{\cdot}\left(\mathrm{d}\MM{q} +\left[\MM{u}\dx{t} +
  \sum_i\bm{\Xi}_i(\bm{x})\circ\mathrm{d} W_i\right]
  \cdot\nabla\MM{q}\right)\dx{^2x} \dx{t},
 \end{split}
\end{align}
where $\ell$ will be the Lagrangian,
$\phi$ is a Lagrange multiplier enforcing the continuity equation, $\MM{q}(\MM{x},t)$ is the back-to-labels map returning the Lagrangian label of the fluid particle at position $\MM{x}$ at time $t$, and $\MM{p}(\MM{x},t)$ is the Lagrange multiplier enforcing the advection of Lagrangian particles. 
\new{The basis functions $\MM{\Xi}_i$ and are tangential to the boundary $\partial\Omega$ of the domain $\Omega$, which we shall assume to be
simply connected.} \\
\\
Following \cite{holm1998hamilton}, we use the specific
Lagrangian for QG
\begin{equation}
\ell(\bm{u},\rho,P )= \int_{\Omega }\left(\frac{1}{2}\rho |\bm{u}|^2-\frac{1}{2}\mathcal{F}\rho \bm{u}\bm{\cdot}\Delta^{-1}\bm{u}+\rho \bm{u}\bm{\cdot}\bm{R}+P (\rho -\rho _0)\right)\dx{^2x}, \label{QG Lagrangian}
\end{equation}
where $\bm{R}$ is the fluid velocity due to the rotation of the planet so that $\hat{\bm{z}}\bm{\cdot}\curl{R}=f$,
and where the  operator $\Delta^{-1}$ is the inverse of the Laplacian operator. \\
\\
After computing the Euler-Lagrange equations and eliminating $\MM{p}$, $\MM{q}$ and $\phi$, computations in  \cite{holm2015variational} lead to the equation
\begin{equation}
\mathrm{d}\left(\frac{1}{\rho _0}\frac{\delta \ell}{\delta \bm{u}}\right)+\mathrm{d}\bm{x}_t\bm{\cdot}\grad{}\left(\frac{1}{\rho _0}\frac{\delta \ell}{\delta \bm{u}}\right)+\frac{1}{\rho _0}\sum_k\frac{\delta \ell}{\delta u^k}\grad{}\mathrm{d}x_t^k+\grad{}\frac{\delta \ell}{\delta \rho }\dx{t}=\bm{0}, 
\ \ \ \frac{\delta \ell}{\delta P}=0,\label{QG EP eqn}
\end{equation}
where $\mathrm{d}\MM{x}_t$ is given in Equation \eqref{eq:dx}
For the QG case this gives
\begin{equation}
  \frac{1}{\rho _0}\frac{\delta \ell}{\delta \bm{u}}=  \bm{u}+ \bm{R}-\mathcal{F}\Delta^{-1}\bm{u}, \quad \frac{\delta \ell}{\delta \rho}
  = \frac{1}{2}|\MM{u}|^2 - \frac{1}{2}\mathcal{F}\MM{u}\cdot\Delta^{-1}\MM{u}
  + \MM{u}\cdot\MM{R} + P, \ \ \ \rho=\rho_0. \label{func deriv wtr velocity}
\end{equation}
Taking the curl of (\ref{QG EP eqn}) and manipulating using vector calculus identities yields
\begin{equation*}
\left(\mathrm{d}+\mathrm{d}\bm{x}_t\bm{\cdot}\grad{}\right)\left[\hat{\bm{z}}\bm{\cdot}\grad{}\times\left(\frac{1}{\rho _0}\frac{\delta \ell}{\delta \bm{u}}\right)\right] = 0. 
\end{equation*}
For the QG reduced Lagrangian we compute
\begin{align*}
\hat{\bm{z}}\bm{\cdot}\grad{}\times\left(\frac{1}{\rho _0}\frac{\delta \ell}{\delta \bm{u}}\right) &= \hat{\bm{z}}\bm{\cdot}(\curl{u})+\hat{\bm{z}}\bm{\cdot}(\curl{R})-\mathcal{F}\Delta^{-1}\hat{\bm{z}}\bm{\cdot}(\curl{u}), \\
& = \lap{\psi }-\mathcal{F}\psi + f,
\end{align*}
after substituting \eqref{func deriv wtr velocity} and introducing the stream function $\psi$ so that $\nabla^\perp\psi=\MM{u}$.
The boundary conditions require that $\psi=0$ on $\partial\Omega$.
The resulting equation is
\begin{equation}
\label{eq:q evo}
  \left(\mathrm{d}+\mathrm{d}\bm{x}_t\bm{\cdot}\grad{}\right)\left[\lap{\psi }-\mathcal{F}\psi +f\right] = 0,
\end{equation}
which is the stochastic QG equation for potential vorticity
$q=\nabla^2\psi - \mathcal{F}\psi + f$. 
This equation has an infinite set of conserved quantities,
\begin{equation*}
  C^p = \int_\Omega q^p \dx{^2x},
\end{equation*}
for $p=1,2,3,\ldots$, with $p=1$ corresponding to the total PV, and $p=2$ is proportional to the enstrophy,  which is given
by $Z=\tfrac{1}{2}\int_\Omega q^2 \dx{^2x}$. 
Although the energy is not conserved, we can deduce that it remains bounded, since
\begin{align*}
  2E & = \int_\Omega \left(|\nabla\psi|^2 + \mathcal{F}\psi^2\right) \dx{^2 x}, \\
  & = \int_\Omega (f-q)\psi\dx{^2 x}, \\
  & \leq \left(\int_\Omega (f-q)^2\dx{^2x}\right)^{1/2}
  \left(\int_\Omega \psi^2\dx{^2 x}\right)^{1/2}, \\
  & \leq C\left(\int_\Omega q^2\dx{^2 x}\right)^{1/2}
  \left(2\int_\Omega |\nabla\psi|^2 + \mathcal{F}\psi^2 \dx{^2x}\right)^{1/2},
\end{align*}
with $C$ a positive constant having used the Poincar\'e inequality, 
and hence
\begin{equation}
  E\leq \sqrt{C/2}\int_\Omega q^2\dx{^2 x},
\end{equation}
and so $E$ is bounded by a constant multiplied by the enstrophy,
a positive constant of motion.

\section{Finite element discretisation}
\label{sec:discrete QG}
In this section we present, for the stochastic QG equation, a discretisation that preserves the total PV and enstrophy.
We use a finite element discretisation, which allows for
the equation to be easily solved on arbitrary meshes, and
in particular on the sphere. \\
\\
The weak form of the stochastic QG equation is obtained by multiplying Equation \eqref{eq:sqg} by a test function $\gamma$, and integrating by parts to obtain
\begin{equation}
  \label{eq:semidiscrete q}
  \mathrm{d} \int_\Omega \gamma q\dx{^2 x}
  - \int_\Omega q\grad{\gamma} \bm{\cdot}
  \left(\grad{}^\perp\psi\dx{t} + \sum_i\bm{\Xi}_i(\bm{x})\circ\mathrm{d} W_i\right)
  \dx{^2 x} = 0,
\end{equation}
where the boundary term vanishes since $\dx{}\MM{x}_t\cdot\MM{n}=0$ on $\partial\Omega$. 
A similar procedure, multiplying the relationship between $\psi$ and $q$ by a test function $\phi$ that vanishes on the boundary leads to
\begin{equation}
\int_\Omega\left( \mathcal{F}\phi \psi +\grad{ \phi}\bm{\cdot} \grad{\psi }\right)\dx{^2x} = \int_\Omega \phi \left(f - q\right)\dx{^2x}. \label{weak form of PV inversion}
\end{equation}
This is just the standard weak form for the Helmholtz equation.\\
\\
We introduce a finite element discretisation by choosing a continuous finite element space $V$, defining
\begin{equation*}
\mathring{V} = \left\{\psi\in V:\psi=0 \mbox{ on
}\partial\Omega\right\}.
\end{equation*}
The finite element discretisation is obtained by choosing $(q,\psi)\in
(V,\mathring{V})$ such that Equations \eqref{eq:semidiscrete
  q}-\eqref{weak form of PV inversion} hold for all test functions
$(\gamma,\phi)\in (V,\mathring{V})$.\\
\\
It follows immediately from the weak form \eqref{eq:semidiscrete q}
that the total PV is conserved, since choosing $\gamma=1$ leads to
\begin{equation*}
  \diff \int_\Omega q \dx{^2x} = 0.
\end{equation*}
The enstrophy is conserved since choosing $\gamma=q$ leads to
\begin{equation*}
  \diff \int_\Omega q^2 \dx{^2x} = 0.
\end{equation*}
\new{Proof that the $p$-th moment is conserved requires 
taking $\gamma=q^{p-1}$, but this is not in $V$ for $p>2$
and thus higher moments are not conserved.}
In the absence of noise, this discretisation also conserves energy; it
reduces to the standard vorticity-stream function finite element formulation.
For time integration we use the implicit midpoint rule, and we obtain
\begin{align*}
\begin{split}
  \int_\Omega \gamma \left(q^{n+1}-q^n\right)\dx{^2 x}   & \\
  \qquad - \int_\Omega \frac{q^{n+1}+q^n}{2} \grad{\gamma} \bm{\cdot}
  \left(\grad{}^\perp\psi^{n+1/2}\Delta t + \sum_i\bm{\Xi}_i(\bm{x})\Delta W_i
  \right)
  \dx{^2 x} & = 0, \quad \forall \gamma \in V, \end{split}\\
  \int_\Omega\left( \mathcal{F}\phi \psi^{n+1/2} +\grad{ \phi}\bm{\cdot} \grad{\psi^{n+1/2} }\right)\dx{^2x} - \int_\Omega \phi \left(f - \frac{q^{n+1}+q^n}{2}\right)\dx{^2x} &=0, \quad \forall \phi \in \mathring{V},
\end{align*}
where $\Delta W_i$ are independent random variables with normal distribution, $N(0,\Delta t)$.
This provides a coupled nonlinear system of equations for
$(q^{n+1},\psi^{n+1/2})$ which may be solved using Newton's method.\\
\\
Taking $\gamma=1$ immediately gives conservation of the total vorticity $\Pi$,
\begin{equation*}
  \int_\Omega \left(q^{n+1}-q^n\right)\dx{^2 x} = 0.
\end{equation*}
Since the implicit midpoint rule conserves all quadratic invariants of the continuous time equations, this scheme conserves the enstrophy $Z$ exactly as well.
The use of the implicit midpoint rule also makes the
scheme unconditionally stable.
The argument behind the bound of the energy of the previous section also holds when restricted to the finite element space, with a constant that is independent of mesh size.

\section{Statistical Properties of the Numerical Scheme}
\label{sec:statmech}
One of our main goals in this work is to understand the properties of the discretisation presented in Section \ref{sec:discrete QG}.
We are particularly motivated by the work of \cite{majda2006nonlinear} and \cite{dubinkina2007statistical},
who looked at the statistical mechanics of discretisations of the deterministic equation but with randomised initial
states.
In particular, \cite{dubinkina2007statistical} looked at how the conservation properties of the discretisation
could affect the statistics.
In both cases, the distribution of states in phase space was given by a Gibbs distribution.
It is therefore of interest whether this same approach could be applied to our discretisation of the stochastic QG equation, and whether this can still be described by the Gibbs distribution.\\
\\
There is a long history of applying statistical mechanics to describe 2D flows, for instance see  \cite{onsager1949statistical}, \cite{kraichnan1967inertial}, \cite{kraichnan1975statistical} or  \cite{kraichnan1980two}.
Some examples of the application to geophysical flows are  \cite{carnevale1981h} and  \cite{carnevale1982statistical},
while there have also been statistical treatments of quasi-geostrophic fluids, such as  \cite{salmon1976equilibrium},  \cite{carnevale1987nonlinear} and  \cite{merryfield2001equilibrium}. 
More recently, the statistical mechanics of numerical discretisations has been considered.
 \cite{majda2006nonlinear} considered Fourier truncations of the QG equation, whilst  \cite{dubinkina2007statistical} considered finite difference methods using Arakawa's Jacobian, conserving energy and enstrophy.
 \cite{dubinkina2007statistical} also considered the other Arakawa schemes that conserve energy but not enstrophy, or conserve enstrophy but not energy, and compared the numerical results with statistics from Gibbs distributions derived under those assumptions.
Since the stochastic QG equation in our finite element discretisation does not conserve energy but does conserve enstrophy, we are actually in exactly this second situation. 
Following that paper, from the maximum entropy principle with constraints of conserved total vorticity $\Pi$ and enstrophy $Z$, we find that the invariant distribution for the finite element discretisation is the Gibbs distribution
$\mathcal{G}(\bm{Q})$, where $\bm{Q}$ is the vector of values describing the discrete $q$ field.  
For our numerical scheme the probability density function for this Gibbs distribution is
\begin{equation}
\mathcal{G}(\bm{Q}) = C^{-1}\exp\left[-\alpha \left(Z(\bm{Q})+\mu\Pi(\bm{Q}) \right)\right], \label{eqn: gibbs distribution}
\end{equation}
where $C$, $\alpha $ and $\mu $ are parameters providing the constraints of conserved $\Pi $, conserved $Z$ and that the integral of the distribution is unity.\\
\\
In this paper, we are interested in computing statistics from this distribution and comparing them to what is obtained from time averages over numerical solutions from the finite element discretisation. 
For example, the expectation of the energy $E$ of the system is then
\begin{equation}
\langle E\rangle = \int_{\mathbb{R}^N} E(\bm{Q})\mathcal{G}(\bm{Q})\dx{\bm{Q}}.
\end{equation}
In general it is not possible to compute this integral analytically.
Our approach is therefore to sample the distribution using a Metropolis algorithm.\\
\\
Before we do this, we will decompose the state vector $\bm{Q}$ into stationary and fluctuating parts:
\begin{equation}
\bm{Q}=\langle \bm{Q}\rangle+\bm{Q}'. \label{eqn: mean and fluct decomp}
\end{equation}
The components $\curlybrackets{Q_i}$ of $\bm{Q}$ are the coefficients in the finite element discretisation, so that for finite element basis $\curlybrackets{\phi_i}$ with $N$ components,
\begin{equation*}
q = \sum_i^NQ_i \phi_i(\bm{x}).
\end{equation*}
Dubinkina and Frank showed in \cite{dubinkina2007statistical} that the  average values $\langle Q_i\rangle$ took a constant value.
Following their computation, we evaluate 
\begin{equation*}
  \left\langle \pfrac{Z}{\bm{Q}}+\mu \pfrac{\Pi }{\bm{Q}}\right\rangle_{\mathcal{G}}
  = \int_{\mathbb{R}^N}\left(\pfrac{Z}{\bm{Q}}+\mu \pfrac{\Pi }{\bm{Q}}\right)C^{-1}\exp\left[-\alpha (Z(\bm{Q})+\mu \Pi (\bm{Q}))\right]\dx{\bm{Q}}.
\end{equation*}
Inspection of the right hand side reveals that 
\begin{equation*}
\left\langle \pfrac{Z}{\bm{Q}}+\mu \pfrac{\Pi }{\bm{Q}}\right\rangle = -\alpha^{-1}\int_{\mathbb{R}^N}\pfrac{}{\bm{Q}}\mathcal{G}(\bm{Q})\dx{\bm{Q}},
\end{equation*}
and if $\mathcal{G}(\bm{Q})$ decays sufficiently fast at infinity then
we conclude that
\begin{equation}
\left\langle \pfrac{Z}{\bm{Q}}+\mu \pfrac{\Pi }{\bm{Q}}\right\rangle = \bm{0}. \label{eqn: av (dZ/dQ + mu *dPi/DQ) = 0}
\end{equation}
In the finite element discretisation with domain $\Omega$, $\Pi(\bm{Q})$ and $Z(\bm{Q})$ are given by
\begin{equation}
\Pi (\bm{Q})  = \sum_i^N\int_\Omega \phi_i(\bm{x}) Q_i\dx{^2x}, \quad Z(\bm{Q})= \frac{1}{2}\sum_{i,j}^N Q_i Q_j \int_\Omega \phi _i(\bm{x}) \phi _j(\bm{x})\dx{^2x}
 \label{discrete_PV and Z}.
\end{equation}
Substituting these into $(\ref{eqn: av (dZ/dQ + mu *dPi/DQ) = 0})$ gives 
\begin{equation*}
\left\langle \int_\Omega \phi_i(\bm{x})(Q(\MM{x})+\mu )\dx{^2x}\right\rangle=0, \hspace{5mm}\forall\phi_i.
\end{equation*}
This means that $Q(\MM{x})$ is the $L^2$-projection of the constant
function $\mu\times 1$ into $V$, but $\mu\times 1\in V$, and we conclude
that $\langle Q(\MM{x)}\rangle = -\mu$ for all $i$. 
For an initial value for the fluid simulation of $\Pi(\bm{Q})=\mathcal{P}_0$, this gives $\langle Q_i\rangle = \mathcal{P}_0/A$, where $A =\int_\Omega \dx{^2x}$.\\
\\
The significance of this result is that it is possible to use the Metropolis algorithm to generate samples from the distribution given in Equation (\ref{eqn: gibbs distribution}) by taking samples from
\begin{equation}
\mathcal{G}'(\bm{Q}') = C^{-1}\exp[-Z(\bm{Q}')].
\end{equation}
This generates samples $\bm{Q}'$ with $\langle \Pi(\bm{Q}')\rangle=0$ and $\langle Z(\bm{Q}')\rangle=\mathcal{Z}'$, which can be transformed into samples of the desired  distribution $\mathcal{G}(\bm{Q})$ (i.e. with $\langle \Pi(\bm{Q})\rangle=\mathcal{P}_0$ and $\langle Z(\bm{Q})\rangle = \mathcal{Z}_0$) by taking
\begin{equation}
Q_i = \frac{\mathcal{P}_0}{A} + Q_i'\sqrt{\frac{\mathcal{Z}_0}{\mathcal{Z}'}-\frac{\mathcal{P}_0^2}{2A\mathcal{Z'}}}. \label{eqn: sample transformation}
\end{equation} 
The Metropolis algorithm can therefore be used to take samples from $\mathcal{G}'(\bm{Q})$, which avoids the evaluation of the parameter $\alpha$ in (\ref{eqn: gibbs distribution}).
The Metropolis algorithm finds to samples of $\mathcal{G}'(\bm{Q}')$ by generating samples from a similar, known distribution, $\mathcal{L}(\bm{Q}')$. 
A given sample  $\bm{Q}'$ is accepted to be from $\mathcal{G}'(\bm{Q}')$ if
 \begin{equation}
 \frac{\mathcal{G}'(\bm{Q}')}{c \mathcal{L}(\bm{Q}')} > 1,
 \end{equation}
 where $c$ is a prescribed constant greater than $1$.
 The Metropolis algorithm also allows us to avoid evaluating the normalisation constants of the distributions.
 A more detailed description of the Metropolis algorithm can be found for example in  \cite{reich2015probabilistic}.\\
 \\
 The known distribution $\mathcal{L}(\bm{Q}')$ that we use is
\begin{equation}
\mathcal{L}(\bm{Q}')=C_\mathcal{L}^{-1}\exp\left[-Z_L(\bm{Q}')\right].
\end{equation}
The lumped enstrophy $Z_L(\bm{Q}')$ is defined as
\begin{equation}
Z_L(\bm{Q}') = \sum_{ij}\int_{\Omega}\frac{1}{2}Q_i^2\phi_i(\bm{x})\phi_j(\bm{x})\dx{^2x}\equiv \sum_i \frac{1}{2} Q_i^2M_i^L,
\end{equation}
with lumped mass matrix $M_i^L:=\sum_j\int_{\Omega}\phi_i(\bm{x})\phi_j(\bm{x})\dx{^2x}$.
This distribution is now straightforward to sample, as
\begin{equation*}
\mathcal{L}(\bm{Q}')= C^{-1}_\mathcal{L}\exp\left[-\sum_i^N \frac{1}{2} Q_i^2M_i^L\right] = C^{-1}_\mathcal{L}\prod_i^N\exp\left[-\frac{1}{2}Q_i^2M_i^L\right].
\end{equation*}
The known distribution is then sampled by generating coefficients $Q'_i$ that are normally distributed with mean 0 and variance $1/M_i^L$. \\
\\
Therefore samples of the Gibbs distribution $\mathcal{G}(\bm{Q})$ are found by using the Metropolis algorithm with the known $\mathcal{L}(\bm{Q}')$ to get samples of $\mathcal{G}'(\bm{Q}')$, and translating the samples $\bm{Q}'$ using (\ref{eqn: sample transformation}).
It is important to note that in general a generated sample
$\bm{Q}$ will not have $\Pi(\bm{Q})=\mathcal{P}_0$ and 
$Z(\bm{Q}) = \mathcal{Z}_0$.
Instead the samples will have $\Pi(\bm{Q})$ and $Z(\bm{Q})$
\textit{distributed} around $\mathcal{P}_0$ and $\mathcal{Z}_0$.
As the resolution increases, the distribution of samples will
become tighter around the possible states in the discretisation.

\section{Numerical Results}
\label{sec:numerics}
In this section we compute numerical trajectories of the finite element discretisation, and compare their time averages with statistics computed using the Gibbs-like distribution. 
The aim was to learn about the properties of the stochastic
QG equation using this approach, and to show that the
Gibbs-like distribution describes the distribution of the
states of discretisation in phase-space, in the limit that
the grid spacing goes to zero.
The code was developed using the Firedrake software suite, which provides code generation from symbolic expressions \cite{firedrake}.
\subsection{Experimental Set-up}
The tests were run on a sphere with unit radius, that is approximated using an icosahedral mesh.
The resolution of this mesh can be refined by subdividing the triangular elements at a given resolution into four smaller triangles to obtain the next resolution.
The refinement level of the icosahedron is the number
of times this process has been repeated to form the mesh.
The function space $V$ was the space of continuous linear
functions on these triangles, whose degrees of freedom (DOFs) are
evaluation of the function are the cell vertices.
The fluid simulation was run for $T$ time steps of $\Delta t = 1$ and
\new{the $\MM{\Xi}$ functions were derived from stream functions give as the
projections of the 
first nine spherical harmonics into the discrete streamfunction space.}
The strength of the stochastic part of the stream function (i.e. the multiplicative constant to the stochastic basis functions) was kept the same for each basis function.
It was found not to affect the average values of the simulation, but increasing it did increase the speed with which the averaged values converged to their limits.\\
\\
For the statistical simulation, the Gibbs-like distribution was used with $\Pi$ and $Z$ corresponding to the initial condition of the fluid simulation.
At the first stage of the statistical simulation, $n$ samples were generated to create an accurate approximation for $\mathcal{Z}'$ to use in equation (\ref{eqn: sample transformation}).
Then the statistical simulation took $n$ samples of the Gibbs-like distribution that were used for the comparison with the statistics of the fluid simulation.
Average properties are found for each simulation by taking the mean value of each sample.
In the case of the fluid simulation, the system at each time step is considered to be a sample.
\subsection{Comparison of Mean Fields}
As the simulations are run, this generates average states in which $q$ is well-mixed over the domain.
Figure 1 shows that \new{the average values predicted from the
fluid simulation and the Monte Carlo simulation lie very close to
one another, even when the fluid simulation starts far
from its average state}. 
This figure plots the average Casimirs  
$\mathcal{C}_3=\int_ \Omega q^3 \dx{^2x}$ and
$\mathcal{C}_4=\int_ \Omega q^4 \dx{^2x}$ as a function of the 
number of samples used in calculating that average for a single 
run of the simulation: we call this the \textit{rolling average}
of the simulations and denote it by angular brackets 
$\langle \cdot \rangle_t$.
This plot was from a run at icosahedron refinement level 4 (corresponding to 2562 $q$ DOFs), 
and the fluid simulation was initialised with $q=\sin\lambda$ for 
latitude $\lambda$.
\begin{figure}[h!]
\centering
\begin{subfigure}{.48\textwidth}
\centering
\includegraphics[width=\textwidth]{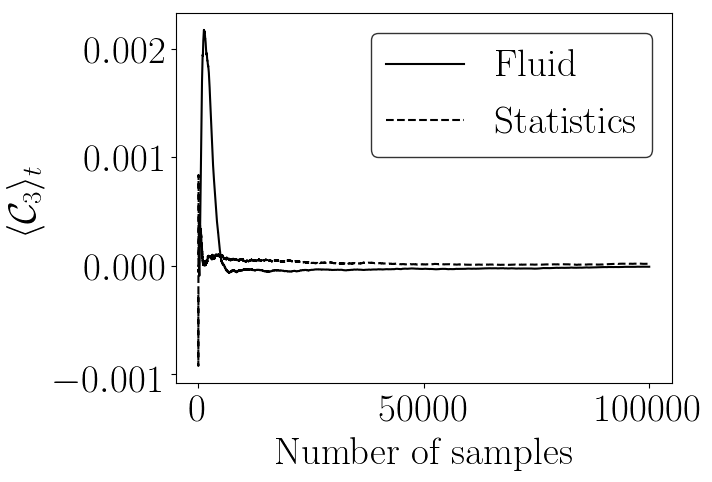} 
\caption*{\small{(a)}}
\end{subfigure}
~~
\begin{subfigure}{.48\textwidth}
\centering
\includegraphics[width=\textwidth]{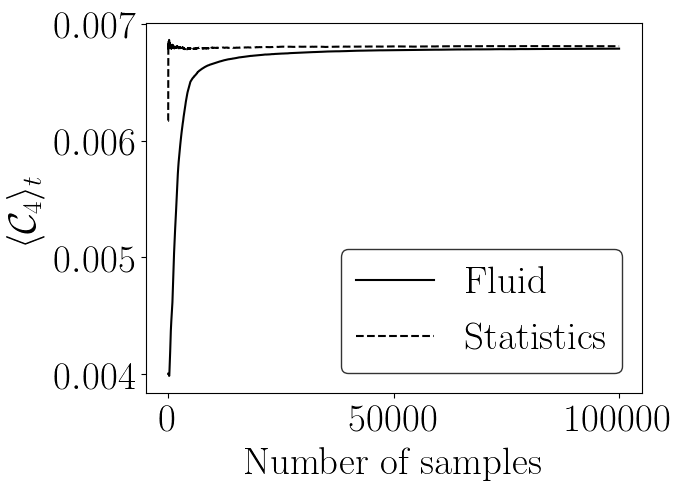}
\caption*{\small{(b)}}
\end{subfigure}
\caption*{\small{\textbf{Figure 1}: The evolution of the rolling
average of the Casimirs (a) 
$\mathcal{C}_3=\int_\Omega q^3 \dx{^2x}$ and (b) 
$\mathcal{C}_4=\int_\Omega q^4\dx{^2x}$ as the stochastic fluid 
simulation is run for $10^5$ time steps and as $10^5$
samples are taken from the Gibbs distribution. 
Both of these are plotted on the same axis, taking one time step of the fluid simulation to be one sample. 
The fluid simulation was run from an initial condition of $q=\sin\lambda$. We observe that both the rolling averages are both converging to the same values, as predicted by the statistical theory.}}
\end{figure}\\
\\
We also took the rolling mean of different diagnostic fields.
These were observed to converge to similar fields for both the 
statistical sampling and from solving the equations of motion.
Figure 2 shows a comparison of the mean $q^2$ fields generated 
via each different method (this is more meaningful than
the mean $q$ field, which shows very little variation over
the domain).
The differences shown \new{over the domain} are due to the 
differences in area of the elements,
whilst the fluid simulation results shown in (a) are 
slightly higher than those in (b) as at finite resolution the 
value of $Z$ of the Gibbs distribution is not exactly equal to
that of the fluid simulation.
As the resolution of the model is refined, the areas of the element will become closer together, and the deviations in the field over the domain should decrease.
This was observed and is plotted in Figure 3.\\
\begin{figure}[h!]
\centering
\begin{subfigure}{.48\textwidth}
\centering
\includegraphics[width=\textwidth]{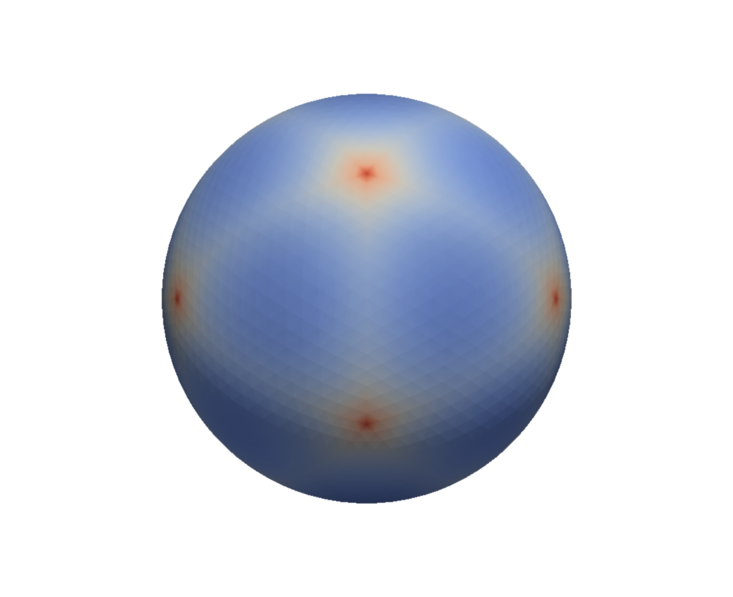} 
\caption*{\small{(a)}}
\end{subfigure}
~~
\begin{subfigure}{.48\textwidth}
\centering
\includegraphics[width=\textwidth]{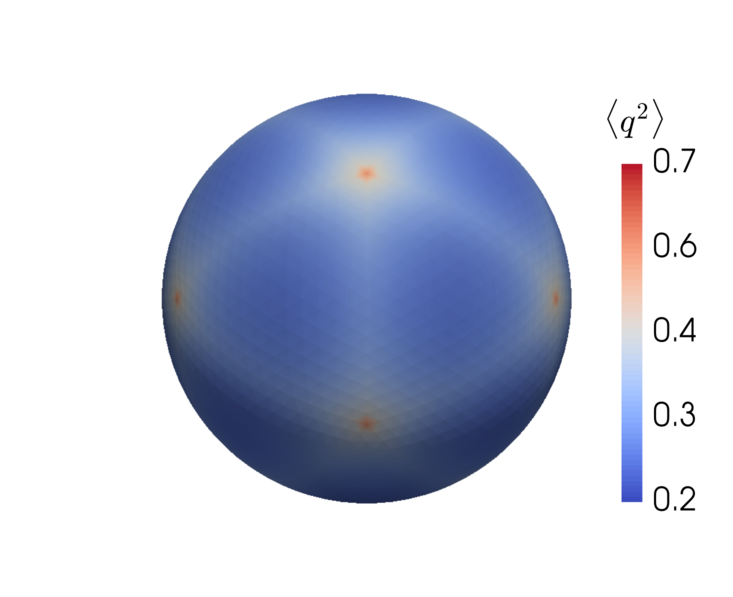}
\caption*{\small{(b)}}
\end{subfigure}
\caption*{\small{\textbf{Figure 2}: A comparison of the mean square potential vorticity $q^2$ field generated at
icosahedron refinement level 4 (
corresponding to 2562 $q$ DOFs) by (a) the stochastic fluid simulation run for 10$^5$ steps with an initially random field and (b) 10$^5$ samples taken from the Metropolis algorithm. We observe that the two fields are essentially the same, as predicted by the statistical theory.}}
\end{figure}
\begin{figure}[h!]
\centering
\begin{subfigure}{.48\textwidth}
\centering
\includegraphics[width=\textwidth]{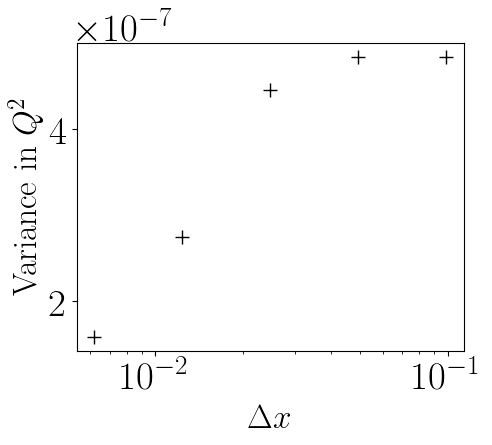} 
\end{subfigure}
~~
\begin{subfigure}{.48\textwidth}
\centering
\includegraphics[width=\textwidth]{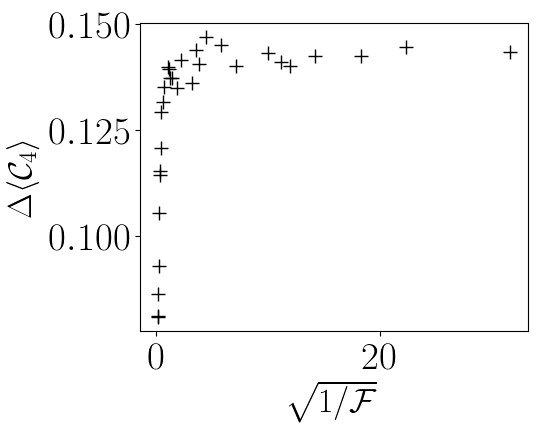}
\end{subfigure}
\caption*{\small{\textbf{Figures 3 and 4}: (Left) A plot of the 
final value of $\int_\Omega (q^2-\bar{q^2})^2\dx{^2x}$
after $10^4$ time steps at $\Delta t=1$ of the fluid simulation, 
where $\bar{q^2}=\int_\Omega q^2\dx{^2x}/A$ and 
$A=\int_{\Omega}\dx{^2x}$. This shows that as the resolution is 
increased, the variations in the mean $q^2$ field reduce.
(Right) The rate of mixing (i.e., rate of convergence of statistics to the equilibrium values) of the fluid simulation as a function of the parameter $\mathcal{F}$, as measured by the speed at which the fourth Casimir $\mathcal{C}_4$ moves away from its initial value. Each point represents the difference at the 2000th time step of the ensemble mean value of $\mathcal{C}_4$ from its initial value. We observe that the rate of mixing is low for large $\mathcal{F}$, which corresponds to large (dimensionless) Rossby deformation radii. This is because scales that are below the Rossby deformation radius are effectively just transported by the flow field without feeding back, and so a large $\mathcal{F}$ places more scales in this category.}}
\end{figure}
\subsection{The Effect of $\mathcal{F}$}
The effect of the constant $\mathcal{F}$ upon the convergence of the model was also investigated.
The equations of motion were solved for a series of different values for $\mathcal{F}$.
For each value of $\mathcal{F}$, the fluid model was run 50 times, creating an ensemble with different realisations of the noise.
Our aim is to learn about the behaviour of the 
discretised stochastic QG equation by initialising the fluid 
at a state far from being well-mixed, and looking at the effect
$\mathcal{F}$ has on the rate of mixing under stochastic noise.
The initial condition was chosen to be $q_0=\sin\lambda$, which is a state far from a well-mixed equilibrium.
This experiment was done with 2562 PV DOFs.
\\
\\
The difference between the initial value of $\mathcal{C}_4$ and its value for the ensemble average was recorded after $T=2000$ time steps.
We found that the average change in $\mathcal{C}_4$ 
over this time was a
good proxy for the rate of mixing of the PV field: 
the larger the difference then the higher rate of mixing.
This was plotted as a function of $1/\sqrt{\mathcal{F}}$, which describes a characteristic length scale.
These values are displayed in Figure 4, which shows smaller differences for smaller length scales (or higher values of $\mathcal{F}$).
This experiment was also performed at different resolutions,
which showed the same behaviour.
\subsection{Convergence With Resolution}
While solving the equations of motion produces only samples with identical $\Pi$ and $Z$, the samples taken from the Gibbs distribution will have different values of $\Pi$ and $Z$, but spread about those values specified in the distribution.
As the resolution is increased, samples taken from the Gibbs distribution should fit more tightly around the samples taken from solving the equations of motion.
\\
\\
To investigate this we ran the fluid simulation and sampled the Gibbs distribution at various resolutions 
to produce histograms of various statistics.
We took 10000 samples from the Gibbs distribution without scaling 
the resultant $\bm{Q}$ field (i.e. the samples had
$\langle\Pi(\bm{Q}\rangle =0 $ and
$\langle Z(\bm{Q}\rangle =Z' $ as described in 
section \ref{sec:statmech}, where $Z'$ will be resolution
dependent).
This will correspond to a fluid simulation with an initial
state of $\Pi=0$ and $Z=Z'$.
Histograms at different resolutions of the $\Pi$ and $Z$ values 
of the statistical samples are plotted in Figure 5.
They have been normalised to remove the resolution dependence
of $Z'$.
The resultant histograms do indeed fit more tightly around
the conserved fluid values as the resolution is increased. 
Figure 6 plots the standard deviation of these histograms as a
function of resolution, which shows that the standard deviation
converges linearly as a function of resolution.
We also would expect to see similar behaviour for other
statistics, although in the case of more complicated statistics
such as variances or higher moments of the distribution, we 
expect that a very large number of samples will be required in 
order to observe this convergence.
\begin{figure}[h!]
\centering
\begin{subfigure}{.48\textwidth}
\centering
\includegraphics[width=\textwidth]{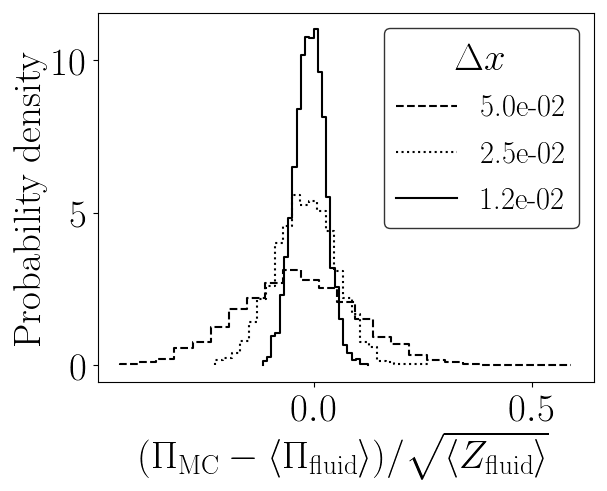}
\caption*{\small{(a)}}
\end{subfigure}
~~
\begin{subfigure}{.48\textwidth}
\centering
\includegraphics[width=\textwidth]{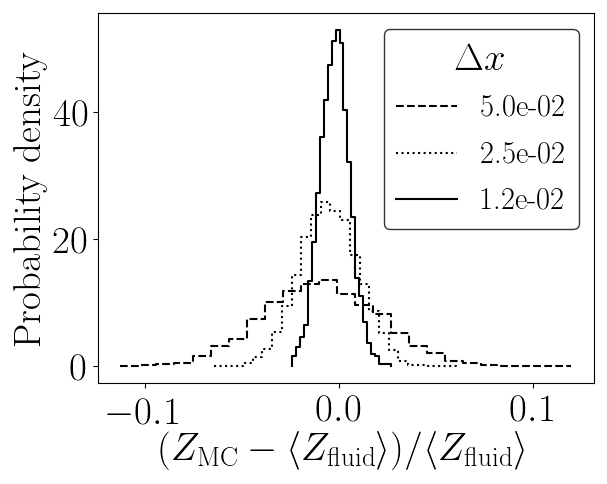} 
\caption*{\small{(b)}}
\end{subfigure}
\caption*{\small{\textbf{Figure 5}: Histograms of the 
(a) $\Pi=\int_\Omega q \dx{^2x}$ and
(b) $Z=\frac{1}{2}\int_\Omega q^2 \dx{^2x}$
values of different samples of the Gibbs distribution at three resolutions.}}
\end{figure}

\begin{figure}[h!]
\centering
\begin{subfigure}{.48\textwidth}
\centering
\includegraphics[width=\textwidth]{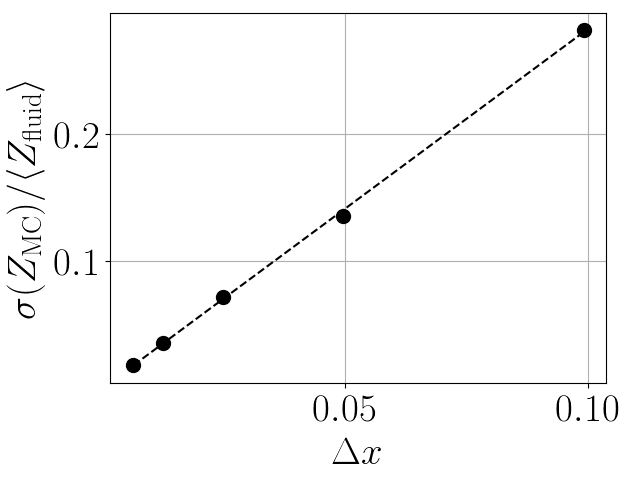}
\caption*{\small{(a)}}
\end{subfigure}
~~
\begin{subfigure}{.48\textwidth}
\centering
\includegraphics[width=\textwidth]{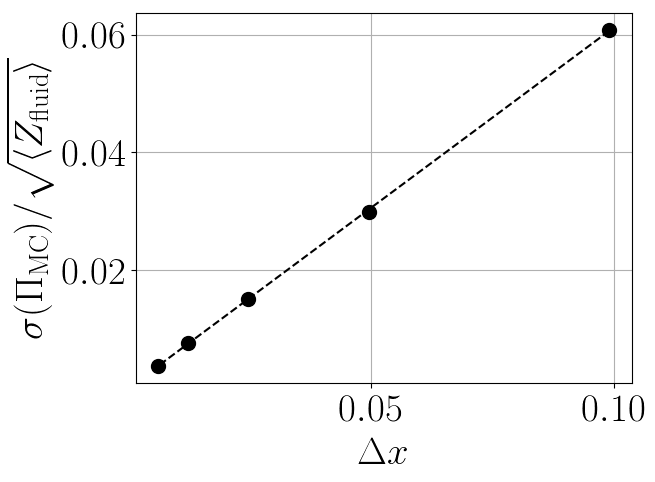} 
\caption*{\small{(b)}}
\end{subfigure}
\caption*{\small{\textbf{Figure 6}: The standard deviations of
the histograms such as those in Figure 5 of
(a) $\Pi=\int_\Omega q \dx{^2x}$ and
(b) $Z=\frac{1}{2}\int_\Omega q^2 \dx{^2x}$ as a function
of resolution. This shows the convergence of the Gibbs distribution to the delta function as the grid spacing
goes to zero.}}
\end{figure}

\subsection{Topography}
The fluid and statistical models also show the same properties when topography is included in the model.
This is done by the addition of an extra term to the potential vorticity definition:
\begin{equation}
q:=\lap{\psi}-\mathcal{F}\psi + f + h.
\end{equation}
We performed similar experiments to those described above by mimicking an isolated mountain, as described as the fifth test case in \cite{williamson1992standard}.
In this case, the topography is described by the following function:
\begin{equation}
h=h_0(1-r/R),
\end{equation}
where $h_0 = 2$, $R=\pi/9$ and $r=\min\left[R^2,(\lambda-\lambda_c)^2+(\theta-\theta_c)^2\right]$, for latitude $\lambda$ and longitude $\theta$.
The centre of the mountain is at $\lambda_c=3\pi/2$ and 
$\theta_c=\pi/6$.
We also completed this experiment with two mountains of the
same radius $R=\pi/9$.
The mountains were both at latitude $\theta_c=\pi/6$, but at
longitudes $\lambda_{c1}=-\pi/4$ and $\lambda_{c2}=\pi/4$.
Figures 7 and 8 show comparisons between the fluid and statistical simulations of mean $\psi$ fields plotted after $T=10^5$ time steps and $n=10^5$ samples were taken.
The third icosahedron refinement level was used, giving 642 $q$ DOFs.
For both the fluid and statistical simulations, the average fields are very similar.
\begin{figure}[h!]
\centering
\begin{subfigure}{.48\textwidth}
\centering
\includegraphics[width=\textwidth]{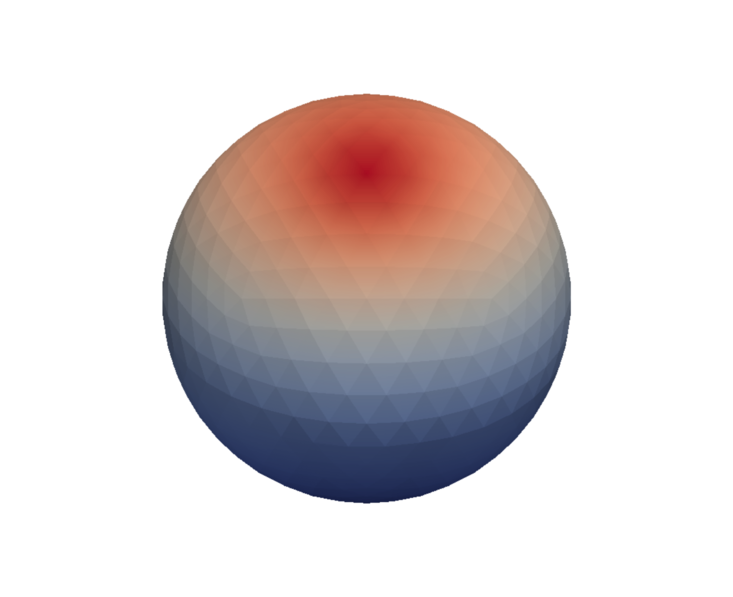} 
\caption*{\small{(a)}}
\end{subfigure}
~~
\begin{subfigure}{.48\textwidth}
\centering
\includegraphics[width=\textwidth]{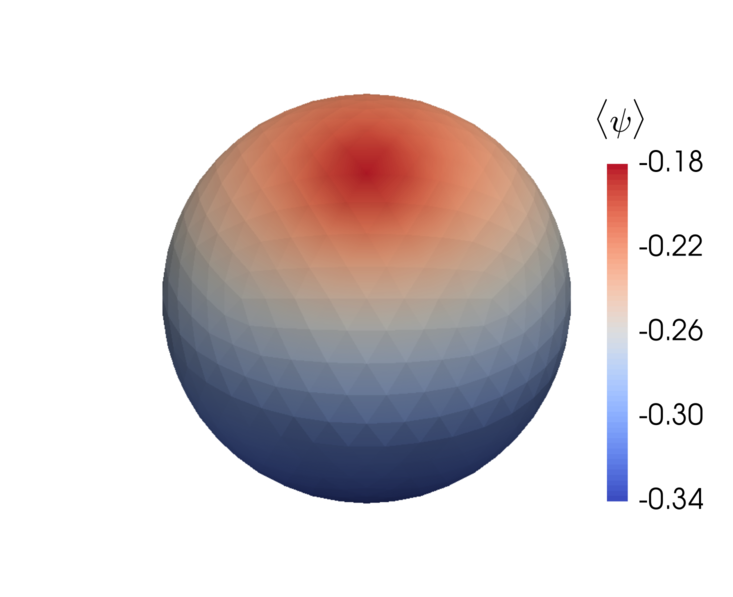}
\caption*{\small{(b)}}
\end{subfigure}
\caption*{\small{\textbf{Figure 7}: A comparison of the mean stream function $\psi$ field generated over topography.
In this case we used flow over a single mountain, and we show
(a) the stochastic fluid simulation run for 10$^6$ steps and (b) 10$^6$ samples taken from the Metropolis algorithm. The plots are essentially the same, illustrating the prediction of the statistical theory.}}
\end{figure}
\begin{figure}[h!]
\centering
\begin{subfigure}{.48\textwidth}
\centering
\includegraphics[width=\textwidth]{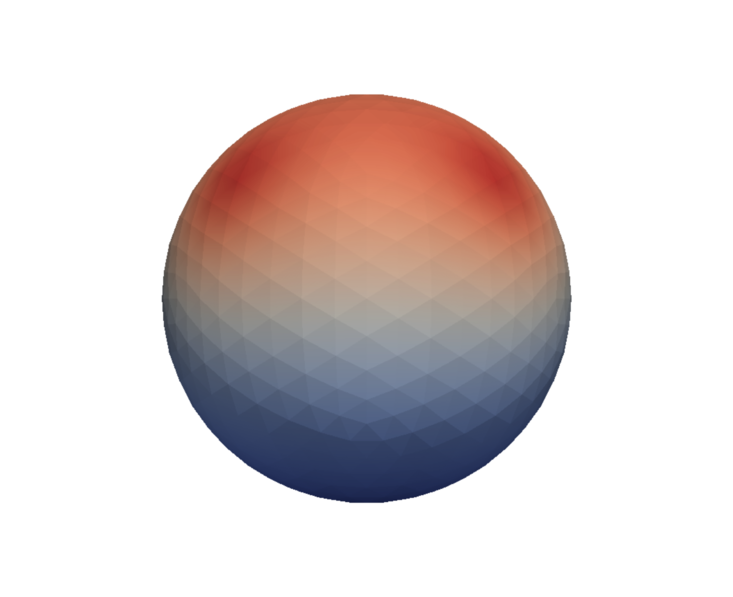} 
\caption*{\small{(a)}}
\end{subfigure}
~~
\begin{subfigure}{.48\textwidth}
\centering
\includegraphics[width=\textwidth]{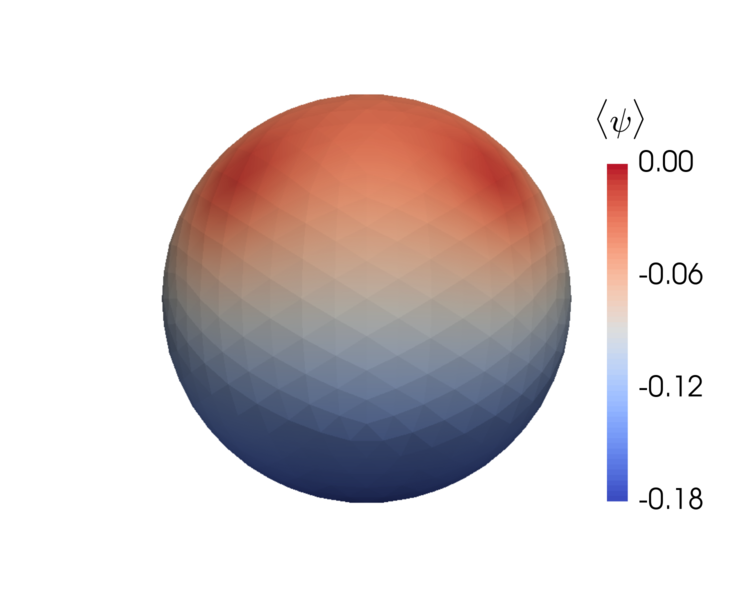}
\caption*{\small{(b)}}
\end{subfigure}
\caption*{\small{\textbf{Figure 8}: A comparison for the two-mountain case of the mean stream function $\psi$ field generated by (a) the stochastic fluid simulation run for 10$^6$ steps and (b) 10$^6$ samples taken from the Metropolis algorithm. The plots are essentially the same, illustrating the prediction of the statistical theory.}}
\end{figure}
\section{Summary and Outlook}
\label{sec:outlook}
Holm's work of \cite{holm2015variational} gave the framework for a stochastic variational principle for QG, which was extended to derive the stochastic QG from the appropriate Lagrangian, before considering quantities conserved by this system. 
A finite element methods was presented to discretise the stochastic QG equation. 
Statistical predictions were then made about this numerical scheme following the approach of \cite{majda2006nonlinear}, which were  tested by sampling using a Metropolis algorithm. 
The statistics generated by sampling the resulting Gibbs-like distribution compared with those averaged quantities from solving the stochastic equations of motion.

\paragraph{Acknowledgements.} TMB was supported by the EPSRC Mathematics of Planet Earth Centre for Doctoral Training at Imperial College London and the University of Reading. CJC was supported by EPSRC grant EP/N023781/1.

\newpage
\bibliography{stochastic-QG} \bibliographystyle{ieeetr}

\end{document}